\documentclass[a4paper,12pt]{article}
\usepackage{amsmath,amsfonts,enumerate,amssymb,amsthm,color}
\usepackage{pgf,tikz}
\title{Arc-transitive cubic abelian bi-Cayley graphs and BCI-graphs}
\author{Hiroki Koike,  \, Istv\'an Kov\'acs\thanks{Corresponding author. \newline
Both authors were supported in part by ARRS - Agencija za raziskovanje Republike Slovenija, program no. P1-0285. \newline 
{\it  E-mail addresses:} hiroki.koike@upr.si (Hiroki Koike),  istvan.kovacs@upr.si (Istv\'an Kov\'acs).}
\\  [+0.75ex]
{\small IAM and FAMNIT, University of Primorska, Muzejski trg 2, 6000 Koper, Slovenia} 
}
\date{}

\newtheorem{thm}{Theorem}[section]
\newtheorem{lem}[thm]{Lemma}

\newtheorem*{thmA}{Theorem A}
\newtheorem*{thmB}{Theorem B}
\theoremstyle{remark}
\newtheorem{rem}[thm]{Remark}
\def\Z{\mathbb{Z}}
\def\S{\mathcal{S}}
\DeclareMathOperator{\aut}{Aut}
\DeclareMathOperator{\bc}{BiCay}
\DeclareMathOperator{\dih}{Dih}
\DeclareMathOperator{\orb}{Orb}
\newcommand{\comment}[1]{}

\begin{document}

\maketitle

\begin{abstract}
A finite simple graph is called a bi-Cayley graph over a group $H$ if it has a semiregular automorphism group, isomorphic to $H,$ 
which has two orbits on the vertex set. 
Cubic vertex-transitive bi-Cayley graphs over abelian groups have been  classified  recently by Feng and 
Zhou (Europ. J. Combin. 36 (2014), 679--693).  In this paper we consider the latter class of graphs and select those in 
the class which are also arc-transitive. 
Furthermore, such a graph is called $0$-type when it is bipartite, and the bipartition classes 
are equal to  the two orbits of the respective semiregular automorphism group. A $0$-type graph can be represented as the graph $\bc(H,S),$ where $S$ is a subset of $H,$ the vertex set of which consists of two copies of $H,$ say $H_0$ and $H_1,$ and the edge set is $\{ \{h_0,g_1\} :  h,g \in H, g h^{-1} \in S\}$.
A bi-Cayley graph $\bc(H,S)$ is called a BCI-graph if for any bi-Cayley graph 
$\bc(H,T),$ $\bc(H,S) \cong \bc(H,T)$ implies that  $T =  h S^\alpha$ for some $h \in H$ and $\alpha \in \aut(H)$. 
It is also shown that every cubic connected arc-transitive $0$-type bi-Cayley graph over an abelian group is a 
BCI-graph.



\end{abstract} 

\section{Introduction}

In this paper all graphs will be simple and finite and all  groups will be finite. 
For a graph $\Gamma$, we let $V(\Gamma)$, $E(\Gamma),$ $A(\Gamma),$ and $\aut(\Gamma)$ denote the 
vertex set, the edge set, the arc set, and the full group of automorphisms of $\Gamma$, respectively.  
A graph $\Gamma$ is called a \emph{bi-Cayley graph} over a group $H$ if it has a semiregular automorphism group, isomorphic to $H,$ which has two orbits in the vertex set. Given such $\Gamma,$ there 
exist subsets $R,L,S$ of $H$ such that $R^{-1}=R,$ $L^{-1}=L,$ $1 \notin R \cup L,$ and $\Gamma \cong \bc(H,R,L,S),$ where the latter graph  is defined to have vertex set the union the 
\emph{right part} $H_0=\{ h_0 : h \in H\}$ and the \emph{left part} $H_1=\{ h_1 : h \in H \};$ and the edge set consists of three 
sets:
\begin{enumerate}[ ]
\item $\{ \{h_0,g_0\} :  g h^{-1} \in R \}$ (\emph{right edges}), 
\item  $\{ \{h_1,g_1\} : g h^{-1} \in L \}$ (\emph{left edges}), 
\item  $\{ \{h_0,g_1\} : g h^{-1} \in S \}$ (\emph{spoke edges}). 
\end{enumerate}
In what follows we will also refer to $\bc(H,R,L,S)$ as a \emph{bi-Cayley representation} of $\Gamma$. 
Regarding bi-Cayley graphs, our notation and terms will follow \cite{FenZ14}.
For the case when $|S| = 1,$ the bi-Cayley graph $\bc(H,R, L, S)$ is also called a \emph{one-matching bi-Cayley graph} 
(see \cite{KovMMM09}).  Also, if $|R| = |L| = s,$ then $\bc(H,R, L, S)$ is said to be an \emph{$s$-type bi-Cayley graph}, 
and if $H$ is abelian, then $\bc(H,R, L, S)$ is simply called an \emph{abelian bi-Cayley graph}. 
If $|L|=|R|=0,$ then $\bc(H,S)$ will be written for $\bc(H,\emptyset,\emptyset,S)$. 
Bi-Cayley graphs have been studied  from various aspects  
\cite{CaoM09,FenZ14,GaoL11,GaoLL11,JinL09,JinL10,JinL11,KovMMM09,LuoG09,WanM06, XuJSL08},  
they have been used by constructions of strongly regular graphs \cite{LeuM93,MalMS07,Mar88,ResJ92} and 
semisymmetric graphs \cite{DuM99a,DuM99b,LuWX06}. 
The cubic vertex-transitive abelian bi-Cayley graphs have been  classified  recently by 
Feng and Zhou \cite{FenZ14} (by a cubic graph we mean a regular graph of valency $3$). %\medskip

In this paper we turn to the class of cubic connected  arc-transitive bi-Cayley graphs over abelian groups. 
Recall that a graph $\Gamma$ is called \emph{arc-transitive} when $\aut(\Gamma)$ is transitive on $A(\Gamma)$. 
From now on  we say that $\Gamma$ is \emph{symmetric} when it is connected and arc-transitive. 
In the first part of our paper we are going to determine the cubic symmetric abelian bi-Cayley graphs.  
Clearly, such a graph is $s$-type for $s \in \{0,1,2\};$ and in fact, the classification  in the case of $0$-type and $2$-type graphs follows from the results in  \cite{ConN07,FenN06,KovMMM09,KutM09}. The respective graphs are listed in Tables 1 and 2.  

\begin{table}
\caption{Cubic symmetric abelian $0$-type bi-Cayley graphs.}
\begin{center}
{\small 
\begin{tabular}{|c|p{2.5in}|p{1in}|c|p{1in}|}   \hline 
no. & $H$                                                                                                        &   $S$                            &   $k$-reg.  & other name \\ \hline                                                                                                        
1. & 
$\Z_{rm} \times \Z_m=\langle a,b \, | \, a^{rm}=b^{rm}=1, b^m=a^{m(u+1)} \rangle,$ 
{\footnotesize $r=3^s p_1^{e_1}\cdots p_t^{e_t},$ $r > 3$  and $ r \ge 11$ if  $m=1,$    
$s \in \{0,1\},$ eve- ry $p_i \equiv 1\pmod 3,$ 
and $u^2+u+1 \equiv 0\pmod r$}           &   $\{1,a,b\}$                &   $1$  & $-$  \\ \hline
2. & 
$\Z_8=\langle a \rangle$                                                                   &   $\{1,a^2,a^3\}$    &   $2$  & 
{\footnotesize M\"obius-Kantor graph} \\ \hline 
3. & 
$\Z_m^2=\langle a,b \rangle,$ $m > 1, m \ne 3$                              &    $\{1,a,b\}$             &   $2$  & $-$ \\ \hline 
4. &  
$\Z_{3m} \times \Z_m=\langle a,b \, | \, a^{3m}=b^{3m}=1, a^m=b^m \rangle,$ 
$m > 1$                                                                                                &   $\{1,a,b\}$   &   $2$  & $-$ \\ \hline 
5. & 
$\Z_3=\langle a \rangle$                                                                   &   $\{1,a,a^{-1}\}$     &   $3$  & $K_{3,3}$ \\ \hline 
6. & 
$\Z_3^2=\langle a,b \rangle$                                                           &   $\{1,a,b\}$             &   $3$  & 
{\footnotesize Pappus graph} \\ \hline 
7. &
$\Z_7=\langle a \rangle$                                                                   &   $\{1,a,a^3\}$             &   $4$  & 
{\footnotesize Heawood graph} \\ \hline 
\end{tabular}
}
\end{center}
\end{table}

\begin{table}
\caption{Cubic symmetric abelian $2$-type bi-Cayley graphs.}
\begin{center}
{\small 
\begin{tabular}{|l|c|c|c|c|p{2in}|}   \hline 
$H$                                              &  $R$              &  $L$                &  $S$        &   $k$-trans  & other name \\ \hline 
$\langle a,b\rangle = \Z_2^2$    &   $\{a,b\}$   & $\{a,b\}$        &  $\{1\}$  &   $2$           & $GP(4,1)$ \\ \hline 
$\langle a \rangle \times \langle b \rangle = \Z_2 \times \Z_{10}$  &  $\{ab^3,ab^{-3}\}$   &    $\{b,b^{-1}\}$     
&    $\{1\}$   &   $2$  &  $-$ \\ \hline 
$\langle a \rangle = \Z_n$           &   $\{a\}$       &   $\{a^k\}$         &  $\{1\}$  &   $2$           & $GP(n,k),$ 
{\footnotesize $(n,k)=(4,1),$ $(8,3),$ $(10,2),$ $(12,5)$, $(24,5)$} \\ \hline
$\langle a \rangle = \Z_n$          &    $\{a\}$      &   $\{a^k\}$     &   $\{1\}$             &  $3$           & $GP(n,k),$ 
{\footnotesize $(n,k)=(5,2),$ $(10,3)$} \\ \hline
\end{tabular}
}
\end{center}
\end{table}

\begin{rem}\label{REM1} 
In order to derive the $0$-type graphs, the key observation is that each such graph is of girth $4$ or $6$. Namely, if $S=\{a,b,c\},$ 
then we find in $\bc(H,S)$ the closed walk: 
$$
( 1_0,a_1,(b^{-1}a)_0,(cb^{-1}a)_1,(b^{-1}c)_0,c_1,1_0 ),
$$ 
here we use that $cb^{-1}a=ab^{-1}c$ holds as $H$ is abelian.  It is a folklore result that the cubic symmetric graphs of girth at most $4$ are $K_4, K_{3,3}$ and $Q_3,$ the graph of the cube. There are infinitely many cubic symmetric graphs of girth $6,$ 
but fortunately, all have been determined in \cite{ConN07,FenN06,KutM09}. A  description of these graphs is given in Theorem \ref{KM}, and using this theorem, it is not hard to deduce Table 1 (see Remark \ref{REM3} for the details).
\end{rem}

\begin{rem}\label{REM2}
The $2$-type bi-Cayley graphs $\bc(\Z_n,\{1,-1\},\{k,-k\},\{0\})$ are also known as the \emph{generalized Petersen graphs}, denoted by $GP(n,k)$.  It was proved by Frucht et al.\ \cite{FruGW71} that $GP(n,k)$ is symmetric for exactly seven pairs: $(4,1),$ $(5,2),$ $(8,3),$ $(10,2), (10,3),$ $(12,5)$ and $(24,5)$. 
Recall that, a  bi-Cayley graph $\bc(H,R,L,S)$ is one-matching if $|S|=1$. 
Symmetric one-matching abelian bi-Cayley graphs are classified in \cite[Theorem 1.1]{KovMMM09}, and 
since the $2$-type cubic bi-Cayley graphs are one-matching,  Table 2 can be read off directly from the latter theorem. 
\end{rem}

In this paper we complete the classification of cubic symmetric abelian bi-Cayley graphs by proving the following theorem:

\begin{thmA}\label{A}
There are exactly four cubic symmetric $1$-type abelian bi-Cayley graphs: $K_4, Q_3, GP(8,3)$ and $GP(12,5)$. 
\end{thmA}

In the second part of this paper we turn to the BCI-property of  cubic symmetric abelian $0$-type bi-Cayley graphs. 
Recall that, these are the graphs in the form $\bc(H,S),$ where $H$ is a finite abelian group and $S$ is a subset of $H$. 
A  bi-Cayley graph $\bc(H,S)$ is said to be a \emph{BCI-graph} if 
for every $\bc(H,T),$ $\bc(H,T) \cong  \bc(H,S)$ implies that $T=h S^\sigma$ for some $h \in H$ and $\sigma \in \aut(H);$  
and the group $H$ is called an \emph{$m$-BCI-group} if every bi-Cayley graph over $G$ of degree at most $m$ is a 
BCI-graph. The study of $m$-BCI-groups was initiated in \cite{XuJSL08}, where it was shown that every group is a $1$-BCI-group, 
and a group is a $2$-BCI-group if and only if it has the property that 
any two elements of the same order are either fused or inverse fused (these groups are described in \cite{LiP97}). 
The problem  of classifying all $3$-BCI-groups  is still open, partial results can be found in 
\cite{JinL09,JinL10,JinL11,KoiK,WieZ07,XuJSL08}. It was proved by the present authors (see \cite[Theorem 1.1]{KoiK})
that the nilpotent $3$-BCI-groups are the groups  $U \times V,$ where $U$ is homocyclic of odd order, and $V$ is trivial, or 
$\Z_{2^r},$ or $\Z_2^r,$ or the quaternion group $\mathbf{Q}_8$  (homocyclic means that it is a direct product of cyclic groups of the same order).  Consequently,  the class of abelian $3$-BCI groups is quite restricted. As our second main result, we prove that 
the situation changes completely when one considers only symmetric graphs. 

\begin{thmB}\label{B}
Every cubic symmetric abelian $0$-type bi-Cayley graph is a BCI-graph. 
\end{thmB}

\section{Preliminaries}

Let $G$ be a group acting on a finite set $V$. For $g \in G$ and $v \in V,$ the image of 
$v$ under $g$ will be written as $v^g$. For a subset $U \subseteq V,$ we will denote by $G_U$ the
point-wise stabilizer of $U$ in $G$, while by $G_{\{U\}}$ the
set-wise stabilizer of $U$ in $G$. If $U = \{u\},$ then $G_u$ 
will be written for $G_{\{u\}}$. 
If $G$ is transitive on $V$ and $\Delta \subseteq V$ is a \emph{block} for $G$ (see \cite[page 12]{DixM96}), then the partition 
$\delta = \{ \Delta^g : g \in G\}$  is called the \emph{system of blocks for $G$ induced by $\Delta$}.
The group $G$ acts on $\delta$ naturally, the corresponding \emph{kernel}  will be denoted by $G_\delta,$ i.e., 
$G_\delta = \{ g \in G :  \Delta'^{\, g}=\Delta' \text{ for all } \Delta' \in \delta \}$. 
For further definitions and results from permutation group theory that will appear later,  we refer the reader to \cite{DixM96}. 
\medskip

Below we collect the main ingredients of this paper. 

\paragraph{2.1 Cubic symmetric graphs.}
For a positive integer $k,$ a \emph{$k$-arc} of a graph $\Gamma$ is an ordered $(k+1)$-tuple $(v_0,v_1,\ldots,v_k)$ of vertices of $\Gamma$ such that, for every $i \in \{1,\ldots,k\},$ $v_{i-1}$ is adjacent to $v_i,$ and for every $i \in \{1,\ldots,k-1\},$ $v_{i-1} \ne v_{i+1}$. The graph $\Gamma$ is called \emph{$(G,k)$-arc-transitive} 
(\emph{$(G,k)$-arc-regular}) if $G$ is transitive (regular) on the set of $k$-arcs of $\Gamma$. If $G = \aut(\Gamma),$ then a $(G,k)$-arc-transitive ($(G,k)$-arc-regular) graph is simply called \emph{$k$-transitive} (\emph{$k$-regular}).
The following result is due to Tutte:

\begin{thm}\label{T}
{\rm \cite{Tut47}} Every cubic symmetric graph is $k$-regular for some $k \le 5$
\end{thm}

In this paper we will occasionally need  information about cubic symmetric graphs of small order, and for this 
purpose use the catalog \cite[Table]{ConD02}. We denote by $FnA, FnB, \ldots$ etc.\ the cubic symmetric graphs on 
$n$ points, and simply write $Fn$ if the graph is uniquely determined by $n$.  
Given an abelian group $G,$ the \emph{generalized dihedral group} $\dih(G)$ is the group $\langle G , \eta \rangle \cong 
G \rtimes \langle \eta \rangle,$ where $\eta$ is an involution and it acts on $G$ as $g^\eta = g^{-1}, g \in G$. 
We have the following description of cubic symmetric graphs of girth $6$: 

\begin{thm}\label{KM}
Let $\Gamma$ be a cubic symmetric graph of girth $6$. Then one of the following holds: 
\begin{enumerate}[(i)]
\item $\Gamma$ is $1$-regular, and $\aut(\Gamma)$ contains a regular normal subgroup isomorphic to 
$\dih(L),$ where $L \cong \Z_{rm} \times \Z_m,  r = 3^s p_1^{e_1} \cdots p_t^{e_t}, r > 3$ and 
$r \ge 11$ if $m=1$, $s \in \{0,1\},$  and every $p_i \equiv 1\pmod 3$.
\item $\Gamma$ is $2$-regular, and $\Gamma \cong GP(8,3),$ or $\aut(\Gamma)$ contains a regular normal subgroup 
isomorphic to $\dih(L),$ where $L \cong \Z_{rm} \times \Z_m,$ $r \in \{1,3\},$ $m > 1,$ and  if $r=1,$ then $m \ne 3$.
\item $\Gamma$ is $3$-regular, and $\Gamma \cong F18$ (the Pappus graph) or $GP(10,3)$ (the Desargues graph).
\item $\Gamma$ is $4$-regular, and $\Gamma \cong F14$ (the Heawood graph).
\end{enumerate}
\end{thm}

\noindent In fact, part (i) is deduced from \cite[Theorem 1.2]{KutM09}, part (ii) from \cite[Theorem 1.1]{KutM09}, and parts 
(iii)-(iv) from \cite[Corollary 6.3]{FenN06} (see also \cite[Theorem 2.3]{ConN07}).

\begin{rem}\label{REM3}
Let $\Gamma$ be a cubic symmetric $0$-type abelian bi-Cayley graph. We are going to show below that 
$\Gamma$ is isomorphic to a graph given in Table 1.  As noted in Remark \ref{REM1}, $\Gamma$ is 
of girth $4$ or $6,$ and it follows that if the girth is equal to $4,$ then $\Gamma \cong K_{3,3}$ or $Q_3$. 
The graph $K_{3,3}$ is isomorphic to the bi-Cayley graph given in row no.\  5 of Table 1, and $Q_3$ 
is isomorphic to the bi-Cayley graph given in row no.\ 3 of Table 1 with $m=2$. 
Assume that  $\Gamma$ is of girth $6$. 
Then by Theorem \ref{KM}, $\Gamma$ is $k$-regular for some $k \le 4$.  We consider each case of the theorem 
separately.  \medskip

\noindent CASE 1. $k=1$. \ In this case $\aut(\Gamma)$ contains a regular normal subgroup $K$ isomorphic to 
$\dih(L),$ where $L \cong \Z_{rm} \times \Z_m,  r = 3^s p_1^{e_1} \cdots p_t^{e_t}, r > 3$ and 
$r \ge 11$ if $m=1,$ $s \in \{0,1\},$  and every $p_i \equiv 1\pmod 3$. 
Consequently, the subgroup $H \le K$ that $H \cong L$ is semiregular and has 
two orbits on $V(\Gamma)$. Notice that, the group $L$ is characteristic in $\dih(L)$. Thus  $H$ is characteristic in $K,$ and 
since $K \trianglelefteq \aut(\Gamma)$ we conclude that $H \trianglelefteq \aut(\Gamma)$.   
Using this and that $\Gamma$ is symmetric, we find that $\Gamma$ is bipartite, and the bipartition classes are 
equal to the orbits of $H$. Therefore, $\Gamma \cong \bc(H,S)$ for a subset $S$ of $H$. 
We may assume without loss of generality that  $1 \in S,$ here $1$ denotes the identity element of $H$. 
 Since $\Gamma$ is arc-transitive and $H$ is normal in $\aut(\Gamma),$ 
there exist $\sigma \in \aut(H)$ and $h \in H$ with the property that $S$ is equal to the orbit of $1$ under 
the mapping $\varphi : x \mapsto x^\sigma h, x \in H$. Thus we may write $S=\{1,a,b\}$ such that $1^\varphi = a,$  
$a^\varphi = b$ and $b^\varphi = 1$. It follows from this that $h=a,$  
$a^\sigma = a^{-1} b$ and $b^\sigma = a^{-1}$.  
This shows that both elements $a$ and $b$ are of the same order. 
On the other hand $\Gamma$ is connected, hence $\langle a,b\rangle = H \cong \Z_{rm} \times \Z_m,$ and 
thus $a$ and $b$ are of order $rm,$ and $\langle a^m \rangle = \langle b^m \rangle$. Then we can write 
$\langle a^m \rangle^\sigma = \langle b^m \rangle^\sigma = 
\langle (b^\sigma)^m \rangle = \langle a^m \rangle,$ and thus $(a^m)^\sigma = (a^m)^u$ for some integer $u,$ $\gcd(u,r)=1$.  
From this $(a^m)^u = (a^m)^\sigma = a^{-m} b^m,$ hence $b^m = a^{m(u+1)}$. 
Also, $(a^m)^{u^2}=(a^m)^{\sigma^2}=(a^{-m}b^m)^\sigma=(a^m)^{-u-1},$ and this  gives that 
$u^2+u+1 \equiv 0 \pmod r$.  To sum up, $\bc(H,\{1,a,b\})$ is one the graphs described in row no.\ 1 of Table 1. 
In fact, any graph in that row is symmetric, the proof of this claim we leave for the reader.  \medskip

\noindent CASE 2. $k=2$. \  In this case $\Gamma \cong GP(8,3),$ or $\aut(\Gamma)$ contains a regular normal subgroup 
isomorphic to $\dih(L),$ where $L \cong \Z_{rm} \times \Z_m,$ $r \in \{1,3\},$ $m > 1,$ and if $r=1,$ then $m \ne 3$. 
We have checked 
by \texttt{Magma} that $GP(8,3)$ admits a bi-Cayley representation given in  row no.\ 2 of Table 1.
Otherwise, copying the same argument as in CASE 1, we derive that $\Gamma \cong \bc(H,S),$ where 
$S = \{1,a,b\},$ and either $H = \langle a,b\rangle \cong \Z_m \times \Z_m,$ $m > 1$ and $m \ne 3,$ or 
$H =\langle a,b \, | \, a^{3m}=b^{3m}=1, a^m=b^m \rangle \cong \Z_{3m} \times \Z_m,$ $m > 1$. 
Therefore, $\bc(H,\{1,a,b\})$ is one of the graphs described in row no.\ 3 of Table 1 in the former case, while 
it is one of the graphs described in row no.\ 4 of Table 1 in the latter case. 
In fact, any graph in these rows is symmetric, the proof is again left for the reader.  \medskip

\noindent CASE 3. $k=3$. \  In this case 
$\Gamma \cong F18$ (the Pappus graph) or $GP(10,3)$ (the Desargues graph). The Pappus graph admits 
a bi-Cayley representation given in  row no.\ 6 of Table 1, and we have checked by \texttt{Magma} that the 
Desargues  graph cannot be represented as a $0$-type abelian bi-Cayley graph.  \medskip

\noindent CASE 4. $k=4$. \  In this case $\Gamma \cong F14$ (the Heawood graph), which admits a bi-Cayley representation 
given in  row no.\ 7 of Table 1. 
\end{rem}

\paragraph{2.2 Quotient graphs.} 
Let $\Gamma$ be an arbitrary finite graph and $G \le \aut(\Gamma)$ which is transitive on $V(\Gamma)$.
For a normal subgroup  $N \triangleleft G$ which is not transitive on $V(\Gamma),$ the \emph{quotient graph} $\Gamma_N$ 
is the graph whose vertices are the $N$-orbits on $V(\Gamma),$ and two $N$-orbits $\Delta_i, i=1,2,$ are adjacent if and only if there exist $v_i \in \Delta_i, i=1,2,$ which are adjacent in $\Gamma$. 
Now, let $\Gamma = \bc(H,R,L,S)$ be a bi-Cayley graph and let $h \in H$. 
Following \cite{FenZ14}, let $R(h)$  denote the permutation of $V(\Gamma)=H_0 \cup H_1$ defined by 
$$
(x_i)^{R(h)} = (x h)_i, \;  x \in H \text{ and } i \in \{0,1\}.
$$ 
We set $R(H) = \{ R(h) : h \in H\}$.  Obviously, $R(H) \le \aut(\Gamma),$ and $R(H)$ is semiregular with orbits $H_0$ and 
$H_1$. Notice that, if $H$ is abelian, then the permutation $\iota$ of $V(\Gamma)$ defined by 
$(x_0)^\iota = (x^{-1})_1$ and $(x_1)^\iota=(x^{-1})_0,  x \in H,$  
is an automorphism of $\Gamma$. Furthermore, the group $\big\langle R(H), \iota \big\rangle  \le \aut(\Gamma)$ is 
regular on $V(\Gamma)$. \medskip

The proof of parts (i)-(ii) of the following lemma can be deduced from \cite[Theorem 9]{Lor84}, and the proof of 
part (iii) is straightforward, hence it is omitted.

\begin{lem}\label{Gamma_N}
Let $\Gamma$ be a cubic symmetric graph and let $N \le \aut(\Gamma)$ be a normal subgroup which has more than 
$2$ orbits on $V(\Gamma)$. Then the following hold:
\begin{enumerate}[(i)]
\item  $\Gamma_N$ is a cubic symmetric graph.
\item $N$ is equal to the kernel of $\aut(\Gamma)$ acting on the system of blocks consisting of the orbits of $N$. 
Moreover, $N$ is regular on each of its orbits. 
\item Suppose, in addition, that $\Gamma = \bc(H,R,L,S)$ where $H$ is an abelian group, and that 
$N = R(K)$  for some  $K <  H$. Then $\Gamma_N$ is isomorphic to the bi-Cayley graph $\bc(H/K,R/K,$ $L/K,S/K)$.
\end{enumerate} 
\end{lem}

\noindent
In part (iii) of the above lemma, $H/K$ denotes the factor group of $H$ by $K$. 
%(notice that, $K$ is normal in $H$ since $R(K)$ is normal in $\aut(\Gamma)$). 
The elements of $H/K$ are the cosets $K h,$ and for a subset $X \subseteq H,$ let 
$X/K$ denote the subset of $H/K$ defined by $X/K = \{ K x : x \in X\}$. 

\paragraph{2.3 Voltage graphs.}  
Let $\Gamma$ be a finite simple graph and $K$ be a finite group whose identity element is denoted by $1_K$. 
For an arc $x = (w,w') \in A(\Gamma)$ we set $x^{-1} = (w',w)$.  A {\em $K$-voltage assignment} 
of $\Gamma$ is a mapping $\zeta : A(\Gamma) \to K$ with
the property $\zeta(x^{-1})=\zeta(x)^{-1}$ for every $x \in A(\Gamma)$.
The values of $\zeta$ are called {\em voltages} and $K$ is called the {\em voltage group}.
Voltages are naturally extended to a directed walk $\vec{W} = (w_1,\ldots,w_n)$ by letting 
$\zeta(\vec{W}) = \prod_{i=1}^{n-1} \zeta((w_i,w_{i+1})).$ 
Fix a  spanning tree $T$ of $\Gamma$. 
Then every edge not in $E(T)$  together with the edges in $E(T)$ span  
a unique circuit of $\Gamma,$ and we shall refer to the circuits obtained in this manner as the {\em base circuits of 
$\Gamma$ relative to $T$}.  The $K$-voltage assignment $\zeta$ is called {\em $T$-reduced} if $\zeta(x) = 1_K$ whenever $x$ is an arc belonging to $A(T)$.  
The {\it voltage graph} $\Gamma \times_{\zeta} K$  is defined to have vertex set $V(\Gamma) \times K,$ and edge set 
\begin{equation}\label{EDGES}
E(\Gamma \times_{\zeta} K) = \Big\{   \{ (w,k),(w',\zeta(x)k) \} :  x=(w,w') \in A(\Gamma) \text{ and } k \in K \Big\}.
\end{equation}
The voltage group $K$ induces an automorphism group of $\Gamma \times_\zeta K$ through the action 
$$
k_{right} : (w,l) \mapsto (w,lk), \; w  \in V(\Gamma) \text{ and } k,l \in K.
$$ 
We set $K_{right} = \{ k_{right} : k \in K\}$. 
Let $g \in \aut(\Gamma \times_\zeta K)$  such that it normalizes 
$K_{right}$. This implies that, if $(w,k) \in V(\Gamma \times_\zeta K)$ and $(w,k)^g = (w',k'),$ then $w'$ does not depend on the choice of $k \in K,$ and the mapping $w \mapsto w'$ is a well-defined permutation of $V(\Gamma)$. The latter permutation 
will be called the {\em projection} of $g,$ obviously, it belongs to $\aut(\Gamma)$.

On the other hand, an automorphism of $\aut(\Gamma)$ is said to {\em lift} to an automorphism of $\Gamma \times_\zeta K$ if it is the projection of some automorphism of $\Gamma \times_\zeta K$. The following ``lifting lemma'' is a special case 
of \cite[Theorem 4.2]{Mal98}:

\begin{thm}\label{LIFT}
Let $\Gamma \times_\zeta K$ be a connected voltage graph, 
where $K$ is an abelian group, and $\zeta$ is a $T$-reduced $K$-voltage assignment. 
Then $\sigma \in \aut(\Gamma)$ lifts to an automorphism of 
$\Gamma \times_\zeta K$ if and only if there exists some $\sigma_* \in \aut(K)$ such that 
for every directed base circuit $\vec{C}$ relative to $T,$
$\sigma_* (\zeta(\vec{C})) = \zeta(\vec{C}^{\, \sigma}).$ 
\end{thm}

For more information on voltage graphs the reader is referred to \cite{GroT87,Mal98}.  

\paragraph{2.4 BCI-graphs.}  For a bi-Cayley graph $\Gamma = \bc(H,S),$ we let $\S(\aut(\Gamma))$ denote the set of all semiregular subgroups of $\aut(\Gamma)$ whose orbits are $H_0$ and $H_1$.  Clearly, $R(H) \in \S(\aut(\Gamma))$ always holds. 
Our main tool in the proof of Theorem B will be the following lemma proved by the present authors:

\begin{lem}\label{B-type}
{\rm \cite[Lemma 2.1]{KoiK}} The following are equivalent for every bi-Cayley graph $\Gamma = \bc(H,S)$:
\begin{enumerate}[(i)]
\item $\bc(H,S)$ is a BCI-graph. 
\item  The normalizer $N_{\aut(\Gamma)}(R(H))$ is transitive on $V(\Gamma),$ and 
every two subgroups in $\S(\aut(\Gamma)),$ isomorphic to $H,$ are conjugate in $\aut(\Gamma)$.
\end{enumerate}
\end{lem}

\section{Proof of Theorem A}

Till the end of the section we keep the following notation: 
$$
\Gamma=\bc(H,\{r\},\{s\},\{1,t\})
$$
 is a cubic symmetric graph, $H = \langle r,s,t \rangle$ is an abelian group, 
and $r,s$ are involutions.  \medskip

The \emph{core} of a subgroup $A$ in a group $B$ is the largest normal subgroup of $B$ contained in $A$. 
In order to derive Theorem A  we analyze the core of $R(H)$ in $\aut(\Gamma)$. 

\begin{lem}\label{L1}
If $R(H)$ has trivial core in $\aut(\Gamma),$ then one of the following holds:
\begin{enumerate}[(i)]
\item $H \cong \Z_2,$ $s=r=t,$ and $\Gamma  \cong K_4$.
\item $H \cong \Z_2^2,$ $s \ne r,$ $t=sr$ and $\Gamma \cong Q_3$.
\end{enumerate}
\end{lem}

\proof
If $\Gamma$ is of girth at most $4,$ then it is isomorphic to $K_4,$ or $K_{3,3},$ or $Q_3$ (see Remark 1.1). 
In the first case we get at once (i), and it is not hard to see that $K_{3,3}$ is impossible. 
Furthermore, we compute by \texttt{Magma} \cite{BosCP97}  that $Q_3$ is possible, $H \cong \Z_2^2,$ 
and $r,s,t$ must be as given in (ii). \medskip

For the rest of the proof we assume that the girth of $\Gamma$ is larger than $4$. 
Then $r \ne s,$ for otherwise, we find the $4$-circuit $(1_0,1_1,r_1,r_0)$.  Then either 
$\langle r,s \rangle \cap \langle t \rangle$ is trivial, and 
\begin{equation}\label{H1}
H = \langle r,s \rangle \times \langle t \rangle \cong \Z_2^2 \times \Z_n;
\end{equation}
or $t$ is of even order, say $2n$, $t^n \in \langle r,s \rangle,$ and 
\begin{equation}\label{H2}
H = \langle r,s,t \rangle \cong \Z_2 \times \Z_{2n}.  
\end{equation} 
Note that, we have $|H|=4n$. 

By Tutte's Theorem (Theorem \ref{T}), $\Gamma$ is $k$-regular for some $k \le 5$.
The order $|\aut(\Gamma)| = |V(\Gamma)| \cdot 3 \cdot 2^{k-1} = 
|H| \cdot 3 \cdot 2^k,$ and thus $|\aut(\Gamma) : R(H)| = 3 \cdot 2^k$.
Consider the action of $\aut(\Gamma)$ on the set of its right $R(H)$-cosets. 
Since $R(H)$ has trivial core in $\aut(\Gamma),$ this action is faithful. Using this and that $R(H)$ acts as a point stabilizer, 
we have an embedding of $R(H)$ into $S_{3 \cdot 2^k-1}$. We shall write below $H \le S_{3 \cdot 2^k-1}$.  
It was proved in \cite[Theorem 1]{BurG89} that, if $n=3m+2$ and $A \le S_n$ is an abelian subgroup, then
\begin{equation}\label{ABEL}
|A| \le  2 \cdot 3^m,
\end{equation}
and equality holds if and only if $A \cong \Z_2 \times \Z_3^m$.  \bigskip 

\noindent CASE 1. $k=1$. \ In this case $\Z_2^2 \le H \le S_5$.  This implies that $|H| =4,$  
$\Gamma \cong Q_3$ (see \cite[Table]{ConD02}), which contradicts that the girth is larger than $4$. \medskip

\noindent CASE 2. $k=2$. \  In this case $H \le S_{11}$. Since $|H|=4n,$ we obtain by \eqref{ABEL} that  $n \le 13$.
We compute by \texttt{Magma} that, if $H$ is given as in \eqref{H1} and $n \le 13$, then 
$\Gamma$ is not edge-transitive.  Furthermore, if $H$ is given as in \eqref{H2} and $n \le 13$, then $\Gamma$ is edge-transitive only if $n=2$ or $n=3$. Consequently,  $\Gamma \cong GP(8,3)$ or $GP(12,5)$ (see \cite[Table]{ConD02}). 
However, we have checked by \texttt{Magma} that in both cases the possible semiregular subgroups have a non-trivial core in the 
full automorphism group, and thus this case is excluded. \medskip

\noindent CASE 3. $k \ge 3$. \   We may assume that $n > 13,$ see the previous paragraph. We find in $\Gamma$ the 
$8$-cycle $( \, 1_0, r_0,  r_1, (rs)_1, (rs)_0, s_0, s_1, 1_1)$.
Thus there must be an $8$-cycle, say $C$, starting with the $3$-arc $(1_0, t_1, t_0, (t^2)_1),$ let this be written in the form:
$$ 
C=\big( \,  1_0, t_1, t_0, (t^2)_1 \, (\delta t^2 )_{x}, (\gamma \delta t^2)_{x'}, (\beta \gamma \delta t^2)_{x''}, 
(\alpha \beta \gamma \delta t^2)_{x'''} \big),
$$
where $x,x',x'',x'''' \in \{0,1\}$ and $\alpha,\beta,\gamma,\delta \in \{1,r,s,t,t^{-1}\}$.
Put $\eta = \alpha \beta \gamma \delta t^2$. Observe that, 
$\eta = t^i  r^j  s^k$ for some integers $i, j, k \ge 0$. Moreover, $i \le 4$ and 
$i=0$ if and only if 
$C = ( \,  1_0, t_1, t_0, (t^2)_1,  (t^2 s)_1,(t s)_0,(t s)_1, s_0 ),$ and so $\eta = s$.
On the other hand, since $1_0 \sim \eta_{x'''}$ and $\eta_{x'''} \ne t_1,$ $\eta \in \{1,r\},$ and we conclude that 
$i > 0$ (recall that $r \ne s$). Now, $1 = \eta^2=t^{2i}r^{2j}s^{2k} = t^{2i},$ which implies that the order of $t$  is at most $8,$ and hence $n \le 8$ 
(see \eqref{H1} and \eqref{H2}), which contradicts that $n > 13$.  This completes the proof of the lemma. \ $\square$

\begin{lem}\label{L2}
Let $R(N)$ be the core of $R(H)$ in $\aut(\Gamma)$. Then one of the following holds:
\begin{enumerate}[(i)]
\item $H = N \times \langle r \rangle,$ and $N r = N s = N t$. 
\item $H = N \times \langle r,s \rangle,$ $r \ne s,$ and $N t = N rs$.
\end{enumerate}
\end{lem}

\proof 
By Lemma \ref{Gamma_N}.(iii), the quotient graph $\Gamma_{R(N)}$ can be written in the form 
$$
\Gamma_{R(N)} = \bc(H/N,\{N r\},\{N s\},\{N,N t\}).
$$ 
We claim that $R(H/N)$ has trivial core in $\aut(\Gamma_{R(N)})$. This and Lemma \ref{L1} will yield (i) and (ii). 

Let $\rho$ be the permutation representation of $\aut(\Gamma)$  derived from its action on the set of 
$R(N)$-orbits.  By Lemma \ref{Gamma_N}.(ii), the kernel $\ker \rho = R(N),$ $\rho(R(H)) = R(H/N),$ and any subgroup of 
$R(H/N)$ is in the form $\rho(R(K))$ for some $N \le K \le H$.  Assume that $\rho(R(K)) \trianglelefteq \aut(\Gamma_{R(N)})$. 
Then $\rho(R(K)) \trianglelefteq \rho(\aut(\Gamma)),$ and hence 
$R(K) \trianglelefteq \aut(\Gamma)$. Thus $R(K) = R(N),$ because $R(N)$ is the core. We find that $\rho(R(K))$ is trivial, 
and the claim is proved. \ $\square$ \medskip

In the next lemma we deal with case (i) of Lemma \ref{L2}. 

\begin{lem}\label{L3}
Let $R(N)$ be the core of $R(H)$ in $\aut(\Gamma),$ and suppose that $N \ne 1$ and case (i) of Lemma \ref{L2} holds. 
Then one of the following holds:
\begin{enumerate}[(i)]
\item $H \cong \Z_2^2,$ $r=s \ne t,$ and $\Gamma \cong Q_3$.
\item $H = \langle r \rangle \times \langle t \rangle \cong \Z_2 \times \Z_4,$ and  
$\Gamma \cong \bc(H,\{r\},\{r t^2\},\{1,t\}) \cong GP(8,3)$.
\end{enumerate}
\end{lem}

\proof In this case $H = N \times \langle r \rangle,$ and $N r = N s = N t$. 
Thus $s= n_1 r,$ and $t=n_2 r$ for some $n_1,n_2 \in N$. Furthermore, $n_1$ is an involution, and since 
$H= \langle r,s,t \rangle,$  $N = \langle n_1,n_2 \rangle$. 

Assume for the moment that $N$ is not a $2$-group, and let $p$ be an odd prime divisor of $|N|$. Then 
$M = \langle n_1,n_2^p \rangle$ is the unique subgroup in $N$ of index $p,$ hence it is  characteristic in $N$.
Using also that $R(N) \trianglelefteq \aut(\Gamma),$ this gives that $R(M) \trianglelefteq \aut(\Gamma)$.
The quotient graph $\Gamma_{R(M)}$  is a cubic symmetric graph on $4 p$ points admitting a $1$-type bi-Cayley 
representation over $H/M$. It was proved in \cite[Theorem 6.2]{FenK07} that $\Gamma_{R(M)}$  is isomorphic to 
one of the graphs: $GP(10,3), GP(10,5),$ and the Coxeter graph $F28$. We compute by \texttt{Magma} 
that none of the these graphs has a $1$-type bi-Cayley representation. We conclude that $N$ is a $2$-group.

Notice that, $N \cong \Z_{2^m}$ or $\Z_2 \times \Z_{2^{m-1}}$. If $|N| \ge 8,$ then $N$ has  a characteristic subgroup 
$M$ such that $|N : M|=8$. Using also that $R(N) \trianglelefteq \aut(\Gamma),$ we find in turn that, 
$R(M) \trianglelefteq \aut(\Gamma),$ and $\Gamma_{R(M)}$ is a cubic symmetric graph on $32$ points which  
admits a $1$-type bi-Cayley representation over $H/M$. Thus $\Gamma$ is isomorphic to 
the Dyck graph $F32$ (see \cite[Table]{ConD02}), which can be excluded by the help of \texttt{Magma}.
Therefore,  $|N| \in \{2,4\},$ and these yield easily cases (i) and (ii) respectively.  \ $\square$ \medskip 

In the next lemma we deal with case (ii) of Lemma \ref{L2}. 

\begin{lem}\label{L4}
Let $R(N)$ be the core of $R(H)$ in $\aut(\Gamma),$ and suppose that $N \ne 1$ and case (ii) of Lemma \ref{L3} holds. 
Then $H = \langle r \rangle \times \langle t \rangle \cong \Z_2 \times \Z_6,$ and 
$\Gamma \cong \bc(H,\{r\},\{r t^3\},$ $\{1,t\}) \cong GP(12,5)$.
\end{lem}

\proof 
In this case  $H = N \times \langle r,s \rangle,$ $r \ne s,$ and $N t = N rs$.
Thus $t =  n_1 r s$ for some $n_1 \in N$.
Since $H = \langle r,s,t \rangle, N = \langle n_1 \rangle$.  
Now, by Lemma \ref{Gamma_N}.(iii) we may write 
$$
\Gamma_{R(N)} = \bc(H/N,\{N r\},\{N s\},\{N, N r s\}) \cong Q_3. 
$$ 

We proceed by defining an $N$-voltage assignment of the quotient graph 
$\Gamma_{R(N)}$. For this purpose we have depicted $\Gamma_{R(N)}$ in Fig.\ 1, where we have also 
fixed the spanning tree $T$ specified by the dashed edges. Now, let $\zeta : A(\Gamma_{R(N)}) \to N$ be 
the $T$-reduced $N$-voltage assignment with its voltages being given in 
Fig.\ 1. To simplify notation we set $\widehat{\Gamma} = \Gamma_{R(N)} \times_\zeta N$.
Recall that $N_{right}$ is a subgroup of $\aut(\widehat{\Gamma})$ (see 2.3). 
Next, we prove the following properties:
\begin{equation}\label{PROP}
\Gamma \cong \widehat{\Gamma}, \text{ and } N_{right}  \trianglelefteq \aut(\widehat{\Gamma}).
\end{equation}

\begin{figure}[t!]
\centering
\begin{tikzpicture}[scale=1.6] ---
%%%%% vertex set
\fill (0,-1) circle (1.5pt); \draw (0,-1) node[below] {\footnotesize $N_0$};
\fill (1,0) circle (1.5pt); \draw (1,0) node[right] {\footnotesize $N_1$};
\fill (0,1) circle (1.5pt); \draw (0,1) node[above] {\footnotesize $(N s)_1$};
\fill (-1,0) circle (1.5pt); \draw (-1,0) node[left] {\footnotesize $(N r)_0$};
\fill (0,0) circle (1.5pt); \draw (0,0) node[below] {\footnotesize $(N rs)_1$};
\fill (1,1) circle (1.5pt); \draw (1,1) node[right] {\footnotesize $(N rs)_0$};
\fill (0,2) circle (1.5pt); \draw (0,2) node[above] {\footnotesize $(N s)_0$}; % $(N r)_0$
\fill (-1,1) circle (1.5pt); \draw (-1,1) node[left] {\footnotesize $(N r)_1$};
%%%%%% spanning tree T
\draw[dashed] (-1,0) -- (0,-1) -- (1,0) -- (0,1) -- (0,2) -- (1,1) -- (0,0) -- (-1,1);
%%%%% voltages 
\draw[->] (0, -1) -- (0,-0.5); \draw (0, -0.5) -- (0,0); \draw (0,-0.5) node[right] {\small $n_1$};
\draw[->] (1, 1) -- (1,0.5); \draw (1,0.5) -- (1,0); \draw (1,0.5) node[right] {\small $n_1$};
\draw (-1,0) -- (-1,1);  \draw (-1,0.5) node[left] {\small $1$};
\draw[->] (-1,0) -- (-0.35,0.65); \draw (-0.35,0.65) --(0,1); \draw (-0.35,0.65) node[above] {\small $n_1$};
\draw[->] (0,2) -- (-0.35,1.65); \draw (-0.35,1.65) --(-1,1); \draw (-0.35,1.65) node[above] {\small $n_1$};

\end{tikzpicture}
\caption{\label{fig2} Voltage assignment $\zeta$ of $\Gamma_{R(N)}$.}
\end{figure}
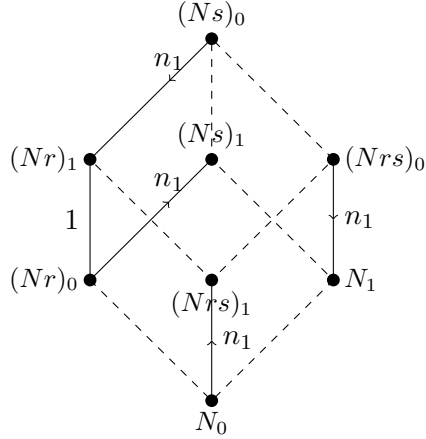

Define the mapping $f : V(\widehat{\Gamma}) \to V(\Gamma)$ by 
$$ 
f \colon ((N x)_0,n) \mapsto (n x)_0 \text{ and } ((N x)_1,n) \mapsto (n x)_1, \; x \in \{1,r,s,rs\}, \; n \in N. 
$$
Notice that, $f$ is well-defined because $\{1,r,s,rs\}$ is a complete set of coset representatives of $N$ in $H$. 
We prove below that $f$ is an isomorphism from $\widehat{\Gamma}$ to $\Gamma$. 
Let $\widehat{v}_1$ and $\widehat{v}_2$  be two adjacent vertices  of $\widehat{\Gamma}$. 
This means that $\widehat{v}_1 = ((N x)_i,n)$  and $\widehat{v}_2= ((N y)_j,\zeta(a) n),$ 
where $a=((N x)_i,( N y)_j)$ is an arc of $\Gamma_{R(N)}$. Then
$ f(\widehat{v}_1) = (x n)_i$ and $f(\widehat{v}_2)= (y \zeta(a)n)_j$.

Let  $i=j=0$. Then it can be seen in Fig.\ 1 that  $y = r x$ and $\zeta(a)=1$. Thus in $\Gamma$ we find 
$f(\widehat{v}_1) = (n x)_0 \sim (r n x)_0 = (y \zeta(a) n)_0 = f(\widehat{v}_2)$. 
Let $i=j=1$. Then $y = s x,$ $\zeta(a)=1,$ and so 
$f(\widehat{v}_1)=(n x)_1 \sim (s n x)_1 = (y \zeta(a) n)_1 = f(\widehat{v}_2)$. Finally, 
let $i=0$ and $j=1$. Then $y = x$ or $y = r s x$. 
In the former case $\zeta(a) = 1,$ and $f(\widehat{v}_1) = (n x)_0 \sim (n x)_1 = (y \zeta(a) n)_1 =  f(\widehat{v}_2)$. 
In the latter case  $\zeta(a) = n_1,$ and 
$$
f(\widehat{v}_1) = (n x)_0 \sim (t n x)_1 = (n_1 r s n x)_1 = (y  \zeta(a) n) = f(\widehat{v}_2).
$$
By these we have proved that $f$ is indeed an isomorphism.
 
For the second part of \eqref{PROP}, compute that $f R(m) f^{-1}$ maps $((N x)_i,n)$ to $((N x)_i,nm)$ for every $m \in N$. 
Thus $f R(m) f^{-1} = m_{righ},$ and so $f R(N) f^{-1} = N_{right}$. Since $R(N) \trianglelefteq \aut(\Gamma),$ 
$N_{right} = f R(N) f^{-1} \trianglelefteq f \aut(\Gamma) f^{-1} = \aut(\widehat{\Gamma}),$ as 
claimed.   \medskip

Now, \eqref{PROP} holds, implying that $\aut(\widehat{\Gamma})$ projects to an edge-transitive 
subgroup of $\aut(\Gamma_{R(N)})$. We obtain from this that the automorphism $\alpha \in \aut(\Gamma_{R(N)})$ lifts, where 
$$
\alpha = (\, (N r)_0, (N r s)_1, N_1)( \, (N r)_1, (N r s)_0, (N s)_1 \, ).
$$
Apply Theorem \ref{LIFT} to $\widehat{\Gamma}$ with $\sigma = \alpha$ and the following  directed base circuits relative to $T$: 
$$ 
\vec{C} =  ( (N s)_0,(N r s)_0,N_1,(N s)_1) \text{ and } \vec{C}' =  (N_0,(N r)_0,(N s)_1,N_1).
$$ 
Let $\sigma_*$ be the automorphism of $N$ given in Theorem \ref{LIFT}. 
Since $\zeta(\vec{C}) = \zeta(\vec{C}')=n_1,$  $\zeta(\vec{C}^{\, \alpha})=\sigma_*(n_1)=\zeta(\vec{C}'^{\, \alpha}),$ which gives $n_1^{-2}=n_1$. Thus $|N|=3$, and this yields easily the statement of the lemma. \ $\square$

\medskip

\noindent {\it Proof of Theorem A.} \  The theorem follows directly from Lemmas \ref{L1} - \ref{L4}. \ $\square$ 

\section{Proof of Theorem B}

Till the end of the section we keep the following notation:
$$ % \begin{equation}\label{GAMMA}
\Gamma=\bc(H,\{1,a,b\})
$$ % \end{equation}
 is a cubic symmetric graph, where $H = \langle a,b \rangle$ is an abelian group.  \medskip 

Recall that, $\S(\aut(\Gamma))$ denotes the set of all semiregular subgroups of $\aut(\Gamma)$ whose orbits are 
$H_0$ and $H_1$. 

\begin{lem}\label{L5}
For every abelian group $X \in \S(\aut(\Gamma)),$  
there exists an involution $\tau_X \in \aut(\Gamma)$ which satisfies the following 
properties:
\begin{enumerate}[(i)] 
\item Every subgroup $Y \le X$ is normalized by $\tau_X$. 
\item The group $\langle X, \tau_X \rangle$ is regular on $V(\Gamma)$. 
\end{enumerate}
\end{lem}

\proof Let $\bc(X,U)$ be a bi-Cayley representation of $\Gamma$ arising from $X$. 
Then as a permutation group of $V(\Gamma)$,  
$X$ is permutation isomorphic to $R(X)$ acting on $V(\bc(X,U))$.  Therefore, it is sufficient to show the existence 
of an involution $\tau \in \aut(\bc(X,U))$ such that  
\begin{itemize}
\item $\tau$ normalizes every $R(Y) \le R(X);$ and
\item the group $\langle R(X), \tau \rangle$ is regular on $V(\bc(X,U))$. 
\end{itemize}

We claim that the permutation $\tau,$ defined by 
$\tau : x_0 \mapsto  (x^{-1})_1$ and $x_1 \mapsto (x^{-1})_0$, satisfies both properties. 
The edge $\{x_0,(u x)_1\}, u \in U,$ is mapped by $\tau$ to the pair $\{ (x^{-1})_1,u^{-1} x^{-1})_0\}$. 
This is an edge of $\bc(X,U),$ and hence $\tau \in \aut(\bc(X,U))$.
Furthermore, a direct computation gives that for $y \in Y,$ $\tau^{-1}  R(y) \tau = R(y^{-1}),$ hence $\tau$ 
normalizes $R(Y),$ and the lemma follows. \ $\square$ 

\begin{lem}\label{L6}
Let $N \le \aut(\Gamma)$ be a normal subgroup such that there exists an $N$-orbit properly contained in $ H_0,$ and let  $X$ be 
an abelian group from $\S(\aut(\Gamma))$. Then  $N <  X$.
\end{lem}

\proof 
Let $\Delta$ be an $N$-orbit such that $\Delta \subset H_0,$ and let us consider $Y = X \cap \aut(\Gamma)_{\{ \Delta \}}$. 
Since $\Delta$ is a block contained in an $X$-orbit, we obtain that $\Delta$ is an $Y$-orbit. We write $\Delta = \orb_Y(v)$. 
Moreover, as $X$ is semiregular, $Y$ is regular on $\Delta$, and by this and Lemma \ref{Gamma_N}.(ii) we have 
\begin{equation}\label{Y=Delta}
|Y| = |\Delta| = |N|.
\end{equation}
 
Let $\tau_X \in \aut(\Gamma)$ be the automorphism defined in Lemma \ref{L5}, and set  $L = \langle X,\tau_X \rangle$. 
According to Lemma \ref{L5} the group $L$  is transitive on $V(\Gamma),$ and also $Y \trianglelefteq L$.  
Denote by $\delta$ the system of blocks  induced by $\Delta$. Then we may write 
$$
\delta = \{  \Delta^l : l \in L \} = \{ \orb_Y(v)^l : l \in L \} = \{ \orb_Y(v^l) : l \in L \}.
$$ 
From this $Y \le \aut(\Gamma)_\delta,$ where $\aut(\Gamma)_\delta$ is the kernel of $\aut(\Gamma)$ acting on $\delta$. 
Since $\aut(\Gamma)_\delta=N,$ see lemma \ref{Gamma_N}.(ii), $Y \le N$. This and \eqref{Y=Delta} imply that  $N = Y < X$. 
\ $\square$ \medskip

For a group $G$ and a prime $p$ dividing $|G|,$ we let $G_p$  denote a Sylow $p$-subgroup of $G$.  
\bigskip

\noindent {\it Proof of Theorem B.} \   We have to show that $\Gamma$ is a BCI-graph.
Let $X \in \S(\aut(\Gamma))$ such that $X \cong H$. 
By Lemma \ref{B-type} and Lemma \ref{L5}, it is sufficient to show the following
\begin{equation}\label{CONJ}
X \text{ and } R(H) \text{ are conjugate in } \aut(\Gamma).
\end{equation}

Recall that the girth of $\Gamma$ is $4$ or $6,$ and if it is $4,$ then $\Gamma$ is isomorphic to 
$K_{3,3}$ or $Q_3$ (see Remark 1.1).  It is easy to see that  \eqref{CONJ} holds when $\Gamma \cong K_{3,3},$ and we 
have checked by the help of \texttt{Magma} that it also holds when  $\Gamma \cong Q_3$. 
Thus assume that $\Gamma$ is of girth $6$.  By Theorem \ref{KM}, $\Gamma$ is $k$-regular for some $k \le 4$.  \medskip

\noindent  CASE 1. $k=1$. \   In this case $\aut(\Gamma)$ contains a regular normal subgroup $K$ isomorphic to 
$\dih(L),$ where $L \cong \Z_{rm} \times \Z_m,  r = 3^s p_1^{e_1} \cdots p_t^{e_t}, r > 3$ and 
$r \ge 11$ if $m=1,$ $s \in \{0,1\},$  and every $p_i \equiv 1\pmod 3$. We have proved in CASE 1 of Remark \ref{REM3}
that $\aut(\Gamma)$ contains a semiregular normal subgroup $N$ such that $N \cong L,$ and the orbits of $N$ are 
$H_0$ and $H_1$.  Notice that, $X$ contains every proper 
characteristic subgroup $K$ of $N$. Indeed, since $N \trianglelefteq \aut(\Gamma),$  
$K \trianglelefteq \aut(\Gamma),$  and Lemma \ref{L6} can be applied for $N,$ implying that $K < X$.  
In particular, if $N$ is not a $p$-group, 
then $N_p <  X$ for every prime $p$ dividing $|N|,$ and thus $N = X$. Since this holds for every $X \in \S(\aut(\Gamma))$ with $X \cong H,$ it holds also for $X=R(H),$ and we get $R(H) = N = X$. In this case \eqref{CONJ} holds trivially. 
Let $N$ be a $p$-group for a prime $p$. Then it follows from the fact that  $N \cong L$ that 
$p > 3,$ and thus both $R(H)$ and $X$ are Sylow $p$-subgroups of $\aut(\Gamma)$. In this case \eqref{CONJ} follows from 
Sylow's Theorem.  \medskip

\noindent  CASE 2. $k=2$. \ In this case $\Gamma \cong GP(8,3),$ or $\aut(\Gamma)$ contains a regular normal subgroup 
isomorphic to $\dih(L),$ where $L \cong \Z_{rm} \times \Z_m,$ $r \in \{1,3\},$ $m > 1,$ and if $r=1,$ then $m \ne 3$. 
If $\Gamma \cong GP(8,3),$ then we have checked by \texttt{Magma} that $H \cong \Z_8$ and \eqref{CONJ} holds. 
Assume that $\Gamma \ncong GP(8,3)$.  We have proved in CASE 2 of Remark \ref{REM3} that $\aut(\Gamma)$ contains a semiregular normal subgroup $N$ such that $N \cong L,$ and the orbits of $N$ are  $H_0$ and $H_1$.
Now, repeating the argument in CASE 1 above, we  obtain that $N=X=R(H)$ if $N$ is not 
a $p$-group.  Let $N$ be a $p$-group for a prime $p$. If $p > 3,$ then both $R(H)$ and $X$ are Sylow $p$-subgroups of 
$\aut(\Gamma),$ and \eqref{CONJ} follows from Sylow's Theorem. We are left with the case that $p \in \{2,3\}$.

Let $p=2$. Since $N \cong L,$ we find that $N \cong \Z_{2^e} \times \Z_{2^e}, e \ge 1$. 
Define $K = \{   x \in N :  o(x) \le 2^{e-1} \}$. Then $K$ is characteristic in $N$ and thus 
$K \trianglelefteq \aut(\Gamma)$. By Lemma \ref{L6}, $K \le X \cap R(H)$. By Lemma \ref{Gamma_N}.(iii),  
the quotient graph $\Gamma_K$ is a $0$-type Bi-Cayley graph over the group $N/K \cong \Z_2^2$.
Then $\Gamma_K \cong Q_3$ and both $N/K$ and $R(H)/K$ are semiregular  on $V(\Gamma_K)$ having orbits the two bipartition 
classes of $\Gamma_K$. Since $X \cong R(H),$  $X/K \cong R(H)/K$. 
A direct computation, using \texttt{Magma}, gives that there are two possibilities: 
$X/K \cong R(H)/K \cong \Z_2^2$ or $\Z_4$. Furthermore, In the former case  
$X/K = R(H)/K,$ which together with $K  <  X \cap R(H)$ yield 
that $X = R(H),$ and \eqref{CONJ} holds trivially. Suppose that the latter case holds and consider 
$\aut(\Gamma)$ acting on the set of $K$-orbits. The kernel of this action is equal to $K,$  see Lemma \ref{Gamma_N}.(ii), 
and thus the image $\aut(\Gamma)/K$ is a subgroup of $\aut(\Gamma_K)$ which is transitive on the set of $2$-arcs of 
$\Gamma_K$. 
However, $\Gamma_K$ is $2$-regular (it is, in fact,  isomorphic to $Q_3$), and we obtain that 
$\aut(\Gamma)/K = \aut(\Gamma_K)$. We compute by \texttt{Magma} that $X/K$ and $R(H)/K$ are conjugate in 
$\aut(\Gamma_K) =  \aut(\Gamma)/K,$ and so \eqref{CONJ} follows from this and the fact that $K < X \cap R(H)$.

Let $p = 3$. Observe first that $|N| > 3$. For otherwise, $\Gamma \cong K_{3,3},$ contradicting that the girth is $6$. 
Since $N \cong L,$ we find that $N \cong \Z_{3^{e+\varepsilon}} \times \Z_{3^e}, e \ge 1,$ $\varepsilon \in \{0,1\},$ and if $\varepsilon = 0,$ then $e \ge 2$. Let $\varepsilon = 0$. 
Define $K = \{   x \in N :  o(x) \le 3^{e-2} \}$. Then $K$ is characteristic in $N$ and thus 
$K \trianglelefteq \aut(\Gamma)$.  By Lemma \ref{L6}, $K \le X \cap R(H)$. By Lemma \ref{Gamma_N}.(iii),  
the quotient graph $\Gamma_K$ is a $0$-type Bi-Cayley graph of the group $N/K \cong \Z_9^2$. 
It follows that $\Gamma_K$ is the unique cubic symmetric graph on $162$ points of girth $6$ 
(see \cite[Table]{ConD02}). 
A direct computation, using \texttt{Magma}, gives that $X/K =  R(H)/K = N/K,$ which together with $K  <  X \cap R(H)$ yield 
that $X = R(H),$ and \eqref{CONJ} holds trivially. 
Let $\varepsilon = 1$.  Define $K = \{   x \in N :  o(x) \le 3^{e-1} \}$.
Then $K$ is characteristic in $N$ and thus 
$K \trianglelefteq \aut(\Gamma)$.  By Lemma \ref{L6}, $K \le X \cap R(H)$. By Lemma \ref{Gamma_N}.(iii),  
the quotient graph $\Gamma_K$ is a $0$-type Bi-Cayley graph of the group $N/K \cong \Z_9 \times \Z_3$. 
It follows that $\Gamma_K$ is the unique cubic symmetric graph on $54$ points (see \cite[Table]{ConD02}). 
A direct computation, using \texttt{Magma}, gives that $X/K =  R(H)/K = N/K,$ which together with $K  <  X \cap R(H)$ yield 
that $X = R(H),$ and \eqref{CONJ} holds also in this case. \medskip

\noindent  CASE 3. $k=3$. \ In this case $\Gamma \cong F18$ (the Pappus graph) or 
$GP(10,3)$ (the  Desargues graph). 
We have checked by \texttt{Magma} that in the former case $H \cong \Z_3^2$ and \eqref{CONJ} holds, and the latter 
case cannot occur. \medskip

\noindent  CASE 4. $k=4$. \  In this case $\Gamma \cong F14$ (the Heawood graph), and 
\eqref{CONJ} follows at once because $X$ and $R(H)$ are Sylow $7$-subgroups of $\aut(\Gamma)$. \ $\square$

\end{document}